\input amstex
\documentstyle{amsppt}
\magnification=\magstep1
\nologo
\TagsOnRight
\font\goth=eusm10
\define\F{\Cal F}
\define\E{\Cal E}
\define\Ii{\hbox{\goth I}}
\define\La{\Cal L}

\define\Oc{\hbox{\goth O}}

\define\Ext{Ext}
\define\a{\alpha}
\define\be{\beta}

\topmatter
\author carlo madonna
\endauthor
\NoRunningHeads
\title a splitting criterion for rank 2 vector bundles on 
hypersurfaces in $\Bbb P^4$ 
\endtitle
\abstract{We show that Horrocks' criterion for the splitting of 
rank two vector bundles in $\Bbb P^3$ can be extended, with some assumptions
on the Chern classes, on non singular hypersurfaces in $\Bbb P^4$.
Extension of other splitting criterions are studied.}\endabstract
\address{Madonna Carlo, 
Pz.Accademia Antiquaria 30, 00147 Roma}\endaddress
\endtopmatter
\email{madonna\@ciop.mat.uniroma3.it}\endemail
\document
{\bf Introduction} \par
A well known criterion of Horrocks (\cite{OSS}, pg.39)
shows that 
a rank 2 vector bundle $\E$ on $\Bbb P^3$
splits as the direct sum of two line bundles 
if and only if $\oplus_n H^1(\Bbb P^3,\E(n))=0$.
It is known that Horrocks criterion fails on a non singular hyperquadric
of $\Bbb P^4$ (see \cite{O}); moreover in \cite{BGS, Remark 3) pg.169}
it is remarked that, for any $r>1$, there exists a non singular
hypersurface $X \subset \Bbb P^4$ of degree $r$
and a non splitting rank two vector bundle $\E$ 
on $X$ with $\oplus_n H^1(\E(n))=0$, for $r>2$, these bundles, form an 
infinite family.

One still can hope that, with a few number of exceptions on the numerical
invariants of $\E$, also on a non singular hypersurface $X \subset \Bbb P^4$
of degree $r$,
the vanishing $\oplus_n H^1(\E(n))=0$ 
forces a rank 2 vector bundle $\E$ on $X$ to split. 

Our main result is the following:

\proclaim{Theorem} 
Let $X \subset \Bbb P^4$ be a non singular hypersurface, of degree $r$;
$\E$ be a rank 2 vector bundle on $X$ with first Chern class $c_1$; let $b$
be the maximum integer such that $h^0(\E(-b)) \ne 0$. Then:
if $\oplus_n H^1(\E(n))=0$, $\E$ splits as the direct sum of two line bundles, unless $-r<2b-c_1<r-2$.
\endproclaim

We remark that the number $-(2b-c_1)$ is a numerical invariant of $\E$,
sometimes called ``degree of stability'': $2b-c_1<0$ is equivalent to 
$\E$ stable (in the sense of Mumford-Takemoto).\par
Examples show that on general hypersurfaces of degree up to $5$, 
one side of the range is ``sharp'' and
``almost sharp'' the other one. Furthermore, we can show that for any $r>5$
there are special smooth hypersurfaces of degree $r$ for which one side of 
the bound is sharp. These bundles are
obtained using an extended version of Serre correspondence 
between subcanonical curves and rank 
two vector bundles.\par
In section $3$ we improve the criterion showing that for
$-r<2b-c_1<r-2$ condition 
$\oplus_n H^1(\E(n))=0$ may be replaced by the vanishing of  
$H^1(\E(n_0))$ for just one suitable $n_0 \in \Bbb Z$, extending in this 
way the Chiantini-Valabrega criterion (see \cite{CV}).\par
We wish to thank Prof.L.Chiantini and Prof.S.Verra for their help in
the preparation of this paper.

\ 

{\bf 0. Notations and Conventions}

With the word "variety" we shall indicate 
an irreducible projective algebraic variety defined
over the complex field $\Bbb C$ while a rank two vector bundle is an 
algebraic rank $2$ vector bundle.
If $Y \subset X$ is a curve,
$\Ii_{Y/X}$, or $\Ii_Y$ if there is no ambiguity, 
denote the ideal sheaf of $Y$ in $X$ and $\omega_Y$ its 
dualizing sheaf (see \cite{H}, ch.III, 7.5) 
If $\F$ is a sheaf on a variety $X$ we define
$h^i(\F):=\dim H^i(X,\F)$. 

When $\text{Pic}(X) \cong \Bbb Z$, we use this isomorphism to identify
line bundles with integers. In particular, for any 
vector bundle $\E$ we consider $c_1(\E)=c_1 \in \Bbb Z$ and write
$\E(n)$ for $\E \otimes \Oc_X(n)$. 
A locally complete intersection curve $Y \subset X$
is called
$a$-subcanonical if $\omega_Y=\Oc_Y(a)$ for some $a \in \Bbb Z$.

We also define 
$$
b(\E)=b=\max \{ n \mid h^0(\E(-n)) \ne 0 \}. 
$$
We say that a rank $2$ vector bundle $\E$ on $X$ ``splits'' if it is
isomorphic to the direct sum
of two line bundles.
We use the notion of stability as given in \cite{OSS,pg. 160}.
When $\text{Pic}(X) \cong \Bbb Z$ this means in our notation that a rank
two vector bundle $\E$ is stable if and only if $2b-c_1<0$.

Note that the number $2b-c_1$ is invariant by twisting
i.e
$$
2b-c_1=2b(\E(n))-c_1(\E(n))
$$ 
$\forall n \in \Bbb Z$. 
If $b=0$, as remarked in \cite{H1, Remark 1.0.1}, 
$\E$ has a global section whose zero-locus has
codimension $2$.

\

{\bf 1.} A well known criterion of Horrocks says that a rank $2$ vector
bundle $\E$ on $\Bbb P^3$
splits if and only if $\oplus_n H^1(\Bbb P^3,\E(n))=0$.
In this section we obtain an extension of this criterion to
rank $2$ vector bundles on a
non singular hypersurface $X \subset \Bbb P^4$. 

\proclaim{Theorem 1.1} 
Let $\E$ be a rank $2$ vector bundle on $X$. Set $r=\text{deg}(X)$. 
If $\oplus_n H^1(\E(n))=0$
then $\E$ splits unless $-r<2b-c_1<r-2$.
\endproclaim

\demo{Proof} As remarked, the number $2b-c_1$ is an invariant of the bundle,
so replacing $\E$ with $\E(-b)$,
we can assume $b=0$. 
We denote with $H$ and $L$ the intersection of $X$ with a generic hyperplane
and a generic plane in $\Bbb P^4$ and with $\E_H,\E_L$
the restrictions of $\E$ to $H$ and $L$. By Bertini theorem there exist
hyperplane sections so that $H$ and $L$ are non singular.

We split proof in two parts. First,  
we suppose that $-c_1 \geq r-2$ and
we prove theorem showing that:
\roster
\item if $H^1(\E(n))=0$ $\forall n \in \Bbb Z$ and $\E_H$ splits, then
$\E$ splits;
\item if $H^1(\E_H(n))=0$ $\forall n \in \Bbb Z$ and $\E_L$ 
splits, then $\E_H$ splits;
\item $\E_L$ splits.
\endroster
When $c_1 \geq r$ we prove something more, showing that
either $h^1(\E(-2))\ne 0$ or $h^1(\E(r-c_1-2))\ne 0$. 

Let us start with $(1)$.
Assume $\E_H=\Oc_H(\a) \oplus \Oc_H(\be)$
with $\a \leq \be$, for example. We have to extend this
isomorphism to $X$.
{}From the exact sequence
$$
0 \to H^0(\E(n-1)) \to H^0(\E(n)) \to H^0(\E_H(n)) \to 0,
$$
the pull back of the sections $s \in H^0(\E_H(-\a))$ and
$s' \in H^0(\E_H(-\be))$
gives two global sections
$t \in H^0(\E(- \a))$ and $t'\in H^0(\E(- \be))$, defining 
an injective morphism
$\mu:\Oc_X(\a) \oplus \Oc_X(\be) \to \E$ 
which is an isomorphism on $H$ since $t,t'$ are indipendent on $H$.
The locus where $\mu$ is not surjective
is the zero-locus of the section $t \wedge t' \in 
H^0(\wedge^2 \E(- \a- \be))$,
which is empty or has codimension one, since
$\wedge^2 \E(- \a- \be)$ is of rank 1.
This locus is necessary empty, otherwise it would intersect
$H$ while the sections $t,t'$ are indipendent on $H$. So $t,t'$ are
indipendent on $X$ and $\E$ splits.

The some argument apply to prove $(2)$. 

To prove $(3)$, note that since we are assuming $b(\E)=0$,
$\E$ has a global section whose zero-locus has codimension two,
say $Y=(s)_0$. Hence $\E_L$ has a nowhere vanishing global section, and this
one define an injective morphism
$s:\Oc_L \to \E_L$, where $\F=\E_L/ \text{Im}(s)$ 
is a subbundle of
$\E_L$. By the exact sequence
$$
0 \to \Oc_L \to \E_L \to \Cal \F \to 0 \tag *
$$
we have
$\text{rk}(\F)=1$ and comparing
Chern classes, we find $\F=\Oc_L(c_1)$. 
Then $(*)$ reads as:
$$
0 \to \Oc_L \to \E_L \to \Oc_L(c_1) \to 0, 
$$
hence $\E_L$ is an extension of $\Oc_L(c_1)$ by
$\Oc_L$ i.e an element of the group
$$
\Ext^1(\Oc_L(c_1),\Oc_L) \cong H^1(\Oc_L(-c_1)).
$$
Using $c_1+r-3<0$ we have 
$h^1(\Oc_L(-c_1))=h^0(\Oc_L(c_1+r-3))=0$
by duality, since $L$ is a plane curve of degree $r$; so 
the only possible extension is the trivial one, i.e
$\E_L= \Oc_L \oplus \Oc_L(c_1)$. 

This proves $(3)$, so one half of the theorem is 
proved.

Assume now that $c_1 \geq r$. In this case $\E$ is stable and so 
it does not split.
Indeed if $\E= \Oc_X(\alpha) \oplus \Oc_X(\beta)$,
being $c_1=\alpha+\beta$, it is $\alpha>0$ or $\beta >0$ and so
$h^0(\E(-1))>0$, absurd since $b=0$.
We will prove that in this case
$h^1(\E(-2)) \ne 0$ or
$h^1(\E(-c_1+r-2)) \ne 0$.
{}From the exact sequence
$$
0 \to \E(-2) \to \E(-1) \to \E_H(-1) \to 0
$$
one gets:
$$
... \to H^0(\E(-1)) \to H^0(\E_H(-1)) \to H^1(\E(-2)) \to ...,
$$
hence if $h^0(\E_H(-1)) \ne 0$ then $h^1(\E(-2)) \ne 0$ too, since
$h^0(\E(-1))=0$. Assume $h^1(\E(-2))=0$.
We will show that $h^1(\E(-c_1+r-2))=0$ implies
$h^0(\E_H(-1)) \ne 0$, proving in this way that if $h^1(\E(-2))=0$ then
$h^1(\E(-c_1+r-2))\ne 0$. 
Indeed assume $h^0(\E_H(-1))=0$; 
{}from 
$$
0 \to \E_H(-2) \to \E_H(-1) \to \E_L(-1) \to 0
$$
we have an exact sequence
$$
... \to H^0(\E_H(-1)) \to H^0(\E_L(-1)) \to H^1(\E_H(-2)) \to ... . 
$$
Suppose we have proved that 
$h^0(\E_L(-1)) \ne 0$. Then $h^1(\E_H(-2)) \ne 0$ 
and the sequence
$$
... \to H^1(\E(-2)) \to H^1(\E_H(-2)) \to H^2(\E(-3)) \to ...
$$
proves that $h^1(\E(-2))\ne 0$ or $h^2(\E(-3)) \ne 0$, absurd since
by Serre duality
$h^2(\E(-3))=h^1(\E(-c_1+r-2))=0$ and by assumption $h^1(\E(-2))=0$.

It remains to show
that $h^0(\E_L(-1)) \ne 0$.
If not, by the exact sequence
$$
0 \to \Oc_L(-1) \to \E_L(-1) \to \Oc_L(c_1-1) \to 0
$$
we have an injective map
$$
H^0(\Oc_L(c_1-1)) \to H^1(\Oc_L(-1)) \cong (H^0(\Oc_L(r-2)))^{\vee}
$$
absurd since
$c_1-1 \geq r-1 \geq 0$.
\qed
\enddemo

By the previous theorem we have the following:

\proclaim{Corollary 1.2} If $\E$ is a non splitting 
rank $2$ vector bundle on $X$ and $2b-c_1 \geq r-2$ then
there is an integer $n_0$ so that $h^1(\E(n_0)) \ne 0$.
\endproclaim

{}From the proof of the theorem, it follows:
\proclaim{Corollary 1.3} Let $\E$ be a rank $2$ vector 
bundle on $X$ and $2b-c_1 \leq -r$. Then
$h^1(\E(-b-2)) \ne 0$ or $h^1(\E(-c_1+b+r-2)) \ne 0$.
\endproclaim

\remark{Remark} If $r=1$ i.e $X=\Bbb P^3$ theorem 1.1 gives Horrocks
splitting criterion.
\endremark

\

{\bf 2.} In this section we will give examples of non splitting rank $2$ vector
bundles with
$\oplus_n H^1(\E(n))=0$ showing that one side of the range of
Theorem 1.1 is ``sharp''. We will use an extended version of Serre
correspondence between curves and bundles.
The same proof of \cite{H1} or \cite{GH} works, 
under the hypothesis $h^1(\La^{\vee})=h^2(\La^{\vee})=0$, which is always
fullfilled when $X$ is a smooth hypersurface in $\Bbb P^4$ and $\La$
is a line bundle on $X$ (for $\text{Pic}(X) \cong \Bbb Z$,
see \cite{H1, Remark 1.1.1}), in this case too. First a definition:

\definition{Definition} Let $\E$ be a rank $2$ vector bundle on a non
singular $3$-dimensional variety, $s \in H^0(\E)$ a global section whose
zero-locus $C=(s)_0$ has codimension $2$. Let $\La \in Pic(X)$ with
$h^1(\La^{\vee})=h^2(\La^{\vee})=0$ and let 
$\phi:\wedge^2 \E \to\La$ be an isomorphism. We say that two triples
$(\E,s,\phi)$ and $(\E',s',\phi')$ are isomorphic if
there is an isomorphism
$\psi:\E \to \Cal E'$ 
and an element
$\lambda \in \Bbb C^*$,
such that
$s'=\lambda \psi (s)$ and $ \phi'=\lambda^2 \phi (\wedge^2 \psi)^{-1}$.
\enddefinition

The correspondence between bundles and curves is given by the following:
\proclaim{Theorem 2.1(Serre)} 
Let $X$ be a non singular 3-dimensional variety.
A curve $C \subset X$
occurs as the zero-locus of a section of a rank $2$ vector bundle 
$\E$ on $X$ if
it is a local complete intersection and
$\omega_C$ is 
isomorphic to the restriction to $C$ of some invertible
sheaf $\La$ on $X$
with
$h^1(\La^{\vee})=h^2(\La^{\vee})=0$.
More precisely, for any fixed invertible sheaf $\La$ on $X$ 
with the previous conditions,
there exist a bijection between the following set of data:\par
$(i)$ the set of triples $(\E,s,\phi)$ modulo isomorphism;\par
$(ii)$ the set of pairs $(C,\xi)$
where $C$ 
is a locally complete intersection curve in $X$, 
and $\xi:\La \otimes \omega_{X} \otimes \Oc_C \to \omega_C$ 
is an isomorphism (modulo $\Bbb C^*$).
\endproclaim
{}From the theorem it follows:
\proclaim{Corollary 2.2} 
Let $C \subset X$ be a locally complete intersection, subcanonical curve and
$\E$ the bundle associated to $C$. Then $C$ is complete intersection
if and only if $\E$ splits.
\endproclaim

\proclaim{Proposition 2.3} 
Let $X$ be a non singular 3-dimensional variety and $\E$ a 
rank $2$ vector bundle on $X$. 
If $\E$ is
generated by global sections, then the zero-locus of a general section 
is non singular.
\endproclaim
\demo{Proof} The same proof of \cite{HM, Theorem 5.1} works in this case too.
\qed
\enddemo

\remark{Remark} 
\roster
\item If $C$ is reduced and connected the bundle $\E$ associated to $C$ 
is unique;
\item if $n \gg 0$ then $\E(n)$ is generated by global sections and 
if $s \in H^0(\E(n))$ is general then
$(s)_0$ is a non singular irreducible curve;
\item if $X \subset \Bbb P^4$ is a non singular 
hypersurface and $\E$ is associated to a curve $C \subset X$
with $\text{deg}(C)=d$ and genus $g$ then
$d=c_2(\E)$, $2g-2=c_2(c_1+r-5)$ and
$C$ is $(c_1+r-5)$-subcanonical.
\endroster
\endremark

\remark{Remark} Let $C \subset \Bbb P^4$ be a curve, complete intersection of
three non singular hypersurfaces of 
degree $l,m,n$, contained in a non singular 
hypersurface $X \subset \Bbb P^4$. Assume that $C$ is
not complete intersection in
$X$ and let $\E$ associated to $C$ on $X$. Then 
$\oplus_n H^1(\E(n))=0$. Indeed using the exact sequences
$$
0 \to {\Ii_{X/ \Bbb P^4}}  \to {\Ii_{C/ \Bbb P^4}}  \to 
{\Ii_{C/X}} \to 0
$$
$$
0 \to \Oc_X \to \E \to {\Ii}_C(c_1) \to 0
$$
we have, 
$$
... \to H^1(\Ii_{C/\Bbb P^4}(n)) \to  H^1(\Ii_{C/X}(n))
\to H^2(\Ii_{X/\Bbb P^4}(n)) \to ...
$$
so that $H^1(\Ii_{C/X}(n))=0$, hence {}from the  
second sequence we have $H^1(\E(n))=0$ for all $n$. 

We call a curve obtained in this way
a "curve of type $(l,m,n)$". 

So if $C \subset \Bbb P^4$ is a curve of type $(l,m,n)$, contained
in a non singular hypersurface $X \subset \Bbb P^4$ with
$\text{deg}(X)=r>n \geq m \geq l>0$, then
$C$ is $(l+m+n-5)$-subcanonical, and if $\E$ is 
associated to
$C$ we have:
\roster
\item $c_1=l+n+m-r$ and $c_2=\text{deg}(C)=lmn$;
\item $\oplus_n H^1(\E(n))=0$.
\endroster
Let $C=(l,m,n)$ as above with $0<l \leq m \leq n <r$. If 
$m+n \leq r$ then $c_1=l+m+n-r \leq l$, since $m+n \leq r$, 
hence $b=0$ and
$2b-c_1=r-l-m-n $.

If $m+n \geq r$ then $c_1=l+m+n-r \geq l$ and $b=m+n-r$, hence
$2b-c_1=m+n-r-l$ since $m-l \geq 0$
and $n-r \geq -\frac{r}2$. 

Observe that in both cases, we always have $2b-c_1 \geq -\frac{r}2$.
\endremark

\example{Example 2.4} Let $X \subset \Bbb P^4$ 
be a general quintic hypersurface,
and let $C \subset X$ be a line, so $C$ is $(-2)$-subcanonical. 
As remarked we have a bundle $\E$ on $X$ associated to $C$ with
$\oplus_n H^1(\E(n))=0$; this shows sharpness on one side , i.e when 
$2b-c_1=2$, for {\bf all quintics}.   
Indeed
by $0 \to H^0(\Oc_X) \to H^0(\E) \to H^0(\Ii_C(-2)) \to 0$
we find $b=0$ hence $2b-c_1=2$. Note that $\E$ does not split
since $C$ is not complete intersection in $X$.
\endexample

\example{Example 2.5} 
Let $C \subset \Bbb P^4$ a line.
By Bertini 
there exist a non singular
hypersurface $X$ containing $C$, of any degree $r>0$.
If $\E$ is associated to $C$,   
as remarked, we have $\oplus_n H^1(\E(n))=0$ and $2b-c_1=r-3$.
So the right bound of Theorem 1.1 is sharp on 
hypersurfaces containing a line.

If $\text{deg}(X)=2$, the range of Theorem 1.1 is sharp. Indeed the only
exception is for $2b-c_1=-1$ and in fact we have
bundles associated to lines, which do not split.
\endexample

\example{Example 2.6} For $2b-c_1$ negative, the best we can do is to consider curves of type $(\frac{r}2,\frac{r}2,\frac{r}2)$ for $r$ even, or
$(\frac{r-1}2,\frac{r-1}2,\frac{r+1}2)$ for $r>1$ odd. {}From the Remark, 
we know that there exists a smooth hypersurface of degree $r$ containing $C$;
the corresponding rank $2$ bundle $\E$ has no cohomology in the middle
and it does not
split. 
As remarked, we have $2b-c_1=r-l-m-n=-\frac{r}2$ for $r$ even or
$2b-c_1=-\frac{r-1}2$ for $r>1$ odd.  

Observe that even if these examples are far {}from our bound yet they
give stable rank $2$ bundles with no cohomology in the middle
(but only on some particular
hypersurface).
\endexample

{\bf 3.} In the present section we weaken the criterion proved above,
showing that condition $\oplus_n H^1(\E(n))=0$ may be replaced
by the vanishing of some $H^1(\E(n_0))$ for just one suitable $n_0$,
extending in this way the Chiantini-Valabrega criterion.
The idea is the
same used in \cite{CV} i.e we find an integer $n_0$
such that:
\roster
\item if $n \geq n_0$ and $H^1(\E(n))=0$ then $H^1(\E(n+1))=0$;
\item if $n \leq n_0$ and $H^1(\E(n))=0$ then $H^1(\E(n-1))=0$.
\endroster

In this section we always twist $\E$ suitably, so that we may
assume $b(\E)=0$. We have to prove a preliminary vanishing:

\proclaim{Proposition 3.1} 
Let $m>r-4$, $H^1(\E(m))=0$ and $H^1(\E(-m-c_1+r-4))=0$.
Then $H^1(\E(m+1))=0$.
\endproclaim

\demo{Proof} By the exact sequence:
$$
0 \to \E(m-1) \to \E(m) \to \E_H(m) \to 0 
$$
and by Serre duality we have the following exact sequence
$$
... \to H^1(\E(m)) \to H^1(\E_H(m)) \to H^2(\E(m-1)) \cong H^1(\E(-m-c_1+r-4))^{\vee} \to ... .
$$  
which, by assumption, gives $h^1(\E_H(m))=0$. 
By the exact sequence: 
$$
... \to H^1(\E(m)) \to H^1(\E(m+1)) \to H^1(\E_H(m+1)) \to ...
$$
using $h^1(\E(m))=0$ it is sufficient to show that
$h^1(\E_H(m+1))=0$.
By 
$$
0 \to \E_H(m) \to \E_H(m+1) \to \E_L(m+1) \to 0
$$
$$
... \to H^1(\E_H(m)) 
\to H^1(\E_H(m+1)) \to H^1(\E_L(m+1)) @>{\alpha}>> H^2(\E_H(m)) \to ...
$$
since $h^1(\E_H(m))=0$,
it is sufficient to show that
$\alpha$ is injective i.e that the dual map 
$$
\alpha^{\vee}:H^0(\E_H(-m-c_1+r-4)) \to H^0(\E_L(-m-c_1+r-4))
$$ 
is
surjective. 
To prove surjectivity of 
$\alpha^{\vee}$ remind that 
$b=0$ and $\text{rk}(\E)=2$, hence
$\E$ has a global section whose zero-locus has codimension two, 
and so for a generic $L$,
$\E_L$ has a nowhere vanishing global section
giving the following exact sequence:
$$
0 \to \Oc_L \to \E_L \to \Oc_L(c_1) \to 0.
$$
{}From this one and by the commutativity of the following diagramm:
$$
\CD
0 @> >> H^0 \Oc_H(-m-c_1+r-4) @> >> H^0 \E_H(-m-c_1+r-4)\\
@. @V{\beta}VV @V{\alpha^{\vee}}VV \\
0 @> >> H^0 \Oc_L(-m-c_1+r-4) @>{\gamma}>> H^0 \E_L(-m-c_1+r-4) 
@> >> 0 \\
@. @V VV \\
 @. 0
\endCD
$$
since $\gamma$ is bijective being $h^0\Oc_L(-m+r-4)=0$ and
$\beta$ is surjective too, since
$L$ is a plane curve hence arithmetically normal, then 
$\alpha^{\vee}$ is surjective.
\qed
\enddemo

Another preliminary vanishing result is the following:
\proclaim{Proposition 3.2} Let $m \leq -c_1$. 
If $H^1(\E(m))=0$ then $H^1(\E(m-1))=0$.
\endproclaim

\demo{Proof} By the exact sequence:
$$
\CD
... \to H^0(\E(m)) @>u>> H^0(\E_H(m)) \to H^1(\E(m-1)) \to ...
\endCD
$$
if $u$ surjects then $h^1(\E(m-1))=0$. 
We assume $u$ not surjective and
we will show that $\E$ splits.
For $n \gg 0$, as remarked, $\E(n)$ has a global section
whose zero-locus is a irreducible non singular curve
$Y$ and
$$
c_2(\E(n))=\text{deg}(C)=c_2+rnc_1+rn^2.
$$
{}From Koszul complex of this section we have an exact sequence
$$
0 \to \Oc_X \to \E(n) \to \Ii_C(c_1+2n) \to 0,
$$
using $h^0(\E(-n))=0$ $\forall n>0$ and $h^0\E \ne 0$,
we find $h^0(\Ii_C(c_1+n)) \ne 0$ hence there is
a non zero sections $s \in H^0(\E)$, giving
a divisor $S$ on $X$ of 
type $(c_1+n)H$ containing $C$, 
and
$c_1+n$ is the minimum integer
$\alpha$ such that
$C \subset \alpha H$. Moreover
$S$ is irreducible, indeed if not
$C$ should be contained in some componenent of $S$, contradicting the
minimality of $c_1+n$.
Let $h$ be a general hyperplane divisor on $H$
and set $\Gamma:=C \cap h \subset S \cap h:=C''$. Note that
$C''$ is irreducible by Bertini theorem since $S$ is irreducible.
Next, if
$h^0(\E_H(m)) \ne 0$,
$\E_H(m)$ has a global section indipendent 
{}from $s$.
By the exact sequence:
$$
0 \to \Oc_H \to \E_H(n) \to {\Ii}_{\Gamma}(c_1+2n) \to 0 \tag *
$$
$\Gamma$ is also contained in a divisor $C'$ on $H$, of type
$(c_1+n+m)h$ indipendent {}from $C'$. 
Then $\Gamma$ is contained
in $C''$ and in $C'$ and so by intersection formula, we find
$$
\deg(\Gamma) \leq (c_1+n)h \cdot (c_1+n+m)h=r(c_1+n)(c_1+n+m)
\leq r(c_1+n)n
$$
since
$m+c_1 \leq 0$. Then
$$
c_2+rnc_1+rn^2 \leq r(c_1+n)n=rnc_1+rn^2
$$
i.e $c_2 \leq 0$. Being $\text{deg}(Y)=c_2 \geq 0$ 
it is
$c_2=0$ and $Y=\varnothing$. 
So $\E(n)$ has a never vanishing global section 
i.e $\E(n)$ splits, proving theorem.
\qed
\enddemo

We can now prove the splitting criterion for unstable rank $2$ vector 
bundles.

\proclaim{Theorem 3.3} Let $\E$ be a non splitting rank $2$ vector bundle and
$c_1 \leq 2-r$. 
\roster
\item If $r+c_1$ is odd then $H^1(\E(\frac{-c_1+r-3}2)) \ne 0$;
\item if $r+c_1$ is even then 
$H^1(\E(\frac{-c_1+r-4}2)) \ne 0$ and $H^1(\E(\frac{-c_1+r-2}2)) \ne 0$.
Moreover if $c_1 \leq -r$ then $H^1(\E(\frac{-c_1+r}2)) \ne 0$ too.
\endroster
\endproclaim
       
\demo{Proof} 
$(2)$ Assume $h^1(\E(\frac{-c_1+r-4}2))=0$. 
By Proposition 3.2 we find
$h^1(\E(n))=0$ for all $n \leq \frac{-c_1+r-4}2$. Setting
$m=\frac{-c_1+r-4}2$ we have
$$
h^1(\E(-m-c_1+r-4))=h^1(\E(\frac{-c_1+r-4}2))=0
$$
and $m>r-4$ so by Proposition 3.1 we have
$h^1(\E(n))=0$ for all $n \geq \frac{-c_1+r-4}2$. It follows
that $h^1\E(n)=0$
$\forall n$ and $\E$ splits by Theorem 1.1, absurd. 
The some argument apply to prove the remaining
non vanishing with the extra hypotesis $c_1 \leq -r$ required for
$h^1(\E(\frac{-c_1+r}2))=0$.

$(1)$ Similar proof works in this case too.
\qed
\enddemo

Now let us turn our attention to the range $c_1 \geq r$. We are able to prove:

\proclaim{Proposition 3.4} Let $m \geq 0$. If $H^1(\E(-m-2))=0$
then $H^1(\E(n))=0$ for all $n \leq -m-2$.
\endproclaim

\demo{Proof}
By the exact sequence
$$
...\to H^0(\E(-m-2)) \to H^0(\E_H(-m-2)) \to H^1(\E(-m-3)) \to 0,
$$
we find $H^1(\E(-m-3)) \cong H^0(\E_H(-m-2))$
since $h^0(\E(-m-2))=0$  
being $-m-2<-1$.
By 
the exact sequence
$$
... \to H^0(\E(-m-1)) \to H^0(\E_H(-m-1)) \to H^1(\E(-m-2)) \to ...,
$$
since
$h^0(\E(-m-1)=h^1(\E(-m-2)=0$ we find $h^0(\E_H(-m-2))=0$, following thesis.
\qed
\enddemo

Appling previous theorem to stable bundles we have the following:

\proclaim{Corollary 3.5} Let $c_1 \geq r$. Then $H^1(\E(-2)) \ne 0$.
\endproclaim
               
\demo{Proof} 
If
$h^1(\E(-2))=0$, 
setting $m=0$ by previous proposition, we find
$h^1(\E(-c_1+r-2))=0$ since $-c_1+r-2 \leq -2$, 
absurd by Corollary 1.3.
\qed
\enddemo

\proclaim{Corollary 3.6} Let $c_1 \geq r$.
\roster
\item If $c_1+r$ is odd and
$H^1(\E(\frac{-c_1+r-3}2))=0$ then
$H^1(\E(n))=0$ for all $n \leq \frac{-c_1+r-3}2$;
\item if $c_1+r$ is even and $H^1(\E(\frac{-c_1+r-4}2))=0$ then
$H^1(\E(n))=0$ for all $n \leq \frac{-c_1+r-4}2$.
\endroster
\endproclaim

\demo{Proof} 
$(1)$
Setting $m=\frac{c_1-r-1}2$ hypothesis reads as
$$
h^1(\E(-m-2))=h^1(\E(\frac{-c_1+r-3}2))=0.
$$
Thesis follows by Proposition 3.4. 

$(2)$ The same argument of $(1)$ applies.
\qed
\enddemo

\remark{Remark} If moreover $c_1 \geq r+1$ 
we can prove in the same way that:\par
if $H^1(\E(\frac{-c_1+r-2}2))=0$ then $H^1(\E(n))=0$ for all $n \leq 
\frac{-c_1+r-2}2$
while if $c_1 \geq r+3$ and
$H^1(\E(\frac{-c_1+r}2))=0$ then $H^1(\E(n))=0$ for all
$n \leq \frac{-c_1+r}2$.
\endremark

\

In the other side, we have:
\proclaim{Proposition 3.7} Let $c_1 \geq r$. 
\roster
\item If $c_1+r$ is odd and $H^1(\E(\frac{-c_1+r-3}2))=0$
then $H^1(\E(n))=0$ for all $n \geq \frac{-c_1+r-3}2$;
\item if $c_1+r$ is even and $H^1(\E(\frac{-c_1+r-4}2))=0$ then
$H^1(\E(n))=0$ for all $n \geq \frac{-c_1+r-4}2$.
\endroster
\endproclaim

\demo{Proof}
$(1)$ 
For $m \gg 0$, as remarked, $\E(m)$ has a global section whose zero-locus is 
an irreducible non singular curve $C$ giving the following exact sequence
$$
0 \to \Oc_X \to \E(m) \to \Ii_C(c_1+2m) \to 0. \tag *
$$
By hypotesis
$$
h^1(\E(\frac{-c_1+r-3}2))=h^1(\Ii_C(\frac{c_1+r-3}2+m))=0,
$$ 
and being $\frac{-c_1+r-3}2<0$ we have
$$
h^0(\E(\frac{-c_1+r-3}2))=h^0(\Ii_C(\frac{c_1+r-3}2+m))=0
$$
too, hence by the exact sequence
$$
0 \to \Ii_C(\frac{c_1+r-3}2+m) \to \Oc_X(\frac{c_1+r-3}2+m ) \to 
\Oc_C(\frac{c_1+r-3}2+m) \to 0
$$
we find
$$
h^0(\Oc_X(\frac{c_1+r-3}2+m))=h^0(\Oc_C(\frac{c_1+r-3}2+m)).
$$
By Riemann-Roch, being $C$ a locally complete intersection and  
$(c_1+2m+r-5)$-subcanonical curve, we see that
$$
h^0(\Oc_C(\frac{c_1+r-3}2+m))-h^0(\Oc_C(\frac{c_1+r-7}2+m))=
d . \tag I
$$
Note that
$$
h^0(\Oc_X(\frac{c_1+r-7}2+m))=h^0(\Oc_C(\frac{c_1+r-7}2+m)),
$$
indeed, $h^0(\E(\frac{-c_1+r-7}2))=0$ since
$\frac{-c_1+r-7}2<-3$ and
by Corollary 3.6, 
we have
$$
h^1(\E(\frac{-c_1+r-7}2))=
h^1(\Ii_C(\frac{c_1+r-7}2+m))=0
$$ 
so previous equality holds and 
(I) reads as:
$$
d=h^0(\Oc_X(\frac{c_1+r-3}2+m))-h^0(\Oc_X(\frac{c_1+r-7}2+m)). \tag II
$$
Now assume $h^1(\E(\frac{-c_1+r-1}2)) \ne 0$. {}From $(*)$ we have
$h^1(\Ii_C(\frac{c_1+r-1}2+m)) \ne 0$ too, hence:
$$
h^0(\Oc_C(\frac{c_1+r-1}2+m))> h^0(\Oc_X(\frac{c_1+r-1}2+m)).
$$
Again using Riemann-Roch we find
$$
h^0(\Oc_C(\frac{c_1+r-1}2+m))-h^0(\Oc_C(\frac{c_1+r-9}2+m))
=2d. \tag III
$$
Note that
$$
h^0(\Oc_X(\frac{c_1+r-9}2+m))=h^0(\Oc_C(\frac{c_1+r-9}2+m)).
$$
Indeed we have $h^0(\Ii_C(\frac{c_1+r-9}2+m))=0$ since 
$\frac{-c_1+r-9}2<\frac{-c_1+r-7}2<0$, and
$h^1(\Ii_C(\frac{c_1+r-9}2+m))=0$ by Corollary 2.6, since
$\frac{-c_1+r-9}2<\frac{c_1+r-3}2$.
Hence (III) is given by:
$$
2d >h^0(\Oc_X(\frac{c_1+r-1}2+m))-h^0(\Oc_X(\frac{c_1+r-9}2+m)). 
\tag IV
$$
Comparing (II) and (IV) we find
$$
\align
2[h^0(\Oc_X &(\frac{c_1+r-3}2+m))-h^0(\Oc_X(\frac{c_1+r-7}2+m))]>\\
& >h^0(\Oc_X(\frac{c_1+r-1}2+m))-h^0(\Oc_X(\frac{c_1+r-9}2+m)).
\tag V
\endalign
$$
By standard calculus we find
$$
\align
h^0(\Oc_X(k))-h^0(\Oc_X(p)) &=\frac1{4!}[4r(k^3-p^3)+ 6r(5-r)(k^2-p^2)+\\
&+2r(2r^2-15r+35)(k-p)]
\endalign
$$
hence (V) gives the following:
$$
2[-\frac{r}{12} +\frac{c_1^2r}4+c_1mr+m^2r+\frac{r^3}{12}]
>\frac{11r}6+\frac{c_1^2r}2+2c_1mr+2m^2r+\frac{r^3}6
$$
and one realizes soon, by standard calculus, that this implies $r<0$, absurd.

$(2)$ The same argument of $(1)$ applies.
\qed
\enddemo

\remark{Remark} 
In the same way we can prove that:\par
if $c_1+r$ is even, $c_1 \geq r$ and $H^1(\E(\frac{-c_1+r-2}2))=0$
then $H^1(\E(n))=0$ \ $\forall n \geq \frac{-c_1+r-2}2$. 
Moreover if
$c_1 \geq r+2$ and
$H^1(\E(\frac{-c_1+r}2))=0$ then $H^1(\E(n))=0$ \ $\forall 
n \geq \frac{-c_1+r}2$.
\endremark

\

{}From Theorem 3.3, Corollary 3.6 and Proposition 3.7 
we have our main splitting criterion:

\proclaim{Theorem 3.8} 
Let $\E$ be a rank $2$ vector bundle on $X$. Then:
\roster
\item if $c_1+r$ is odd, $\E$ splits if and only if
$H^1(\E(\frac{-c_1+r-3}2))=0$;
\item if $c_1+r$ is even, 
$\E$ splits if and only if
$H^1(\E(\frac{-c_1+r-4}2))=0$,
\endroster
unless $-r<-c_1<r-2$.
\endproclaim

\remark{Remark} 
If $c_1 \leq 2-r$ or $c_1 \geq r+2$ 
we find another suitable vanishing: $\E$ splits if and only if
$H^1(\E(\frac{-c_1+r-2}2))=0$. 
Moreover if
$c_1 \leq -r$ or $c_1 \geq r+4$ then $\E$ splits if and only if 
$H^1(\E(\frac{-c_1+r}2))=0$.

Note that for $r=1$ i.e if $X=\Bbb P^3$ we obtain the Chiantini-Valabrega
splitting criterion (\cite{CV}), 
while for $r=2$ we get a result of Ottaviani (\cite{O, Theorem 3.8}).
\endremark

\remark{Remark} Results presented in this paper concern obviously, via Serre
correspondence, complete intersection and projective normality. Indeed 
by the exact sequences
$$
0 \to \Oc_X \to \E \to \Ii_C(c_1) \to 0
$$        
if $C$ is $a$-subcanonical then
$\omega_C=\Oc_C(a)=\Oc_C(c_1+r-5)$
so previous theorem reads as follow: 
$C$ is complete intersection if and only if:
\roster
\item $H^1(\Ii_C(\frac{a}2+1))=0$, if $a$ is even; 
\item $H^1(\Ii_C(\frac{a+1}2))=0$, if $a$ is odd,
\endroster
unless $-3<a<2r-5$.
\endremark

\

The following results are a natural generalization of the previous ones and
can be proved by a straightforward arrangement of the previous procedure.

\proclaim{Theorem 3.9} Let $X$ be a non singular 3-dimensional variety. 
Assume $\text{Pic}(X) \cong \Bbb Z$,$h^1(\Oc_X(n))=0$ for all $n \in \Bbb Z$.
If $\omega_X= \Oc_X(e)$ then: if
$\oplus_n H^1(\E(n))=0$ then $\E$ splits unless $-e-5<2b-c_1<e+3$.
\endproclaim

\proclaim{Theorem 3.10} Let $X$ be as above and $\E$ a rank $2$ vector 
bundle on $X$. Then
$\E$ splits if and only if:
\roster
\item $H^1(\E(\frac{-c_1+e+2}2))=0$, if
$c_1+e$ is even;
\item $H^1(\E(\frac{-c_1+e+1}2))=0$, if $c_1+e$ is odd.
\endroster
unless $-e-5<2b-c_1<e+3$.
\endproclaim


\widestnumber\key{FLgo}
\Refs\nofrills{REFERENCES}

\ref \key CV \by L.Chiantini-P.Valabrega \paper Subcanonical 
curves and complete intersections
in projective 3-Space \jour Ann. Mat. Pura Appl. 
(4) 
\vol 138 \yr 1984 \pages {309--330}  \endref

\ref \key BGS \by R.O.Buchweitz-G.M.Greuel-F.O.Schreyer
\paper Cohen-Macaulay modules on hypersurface singularities II
\jour Invent. Math.
\vol 88 \yr 1987 \pages{165--182}
\endref

\ref \key {GH} \by P.Griffiths-J.Harris
\paper Two proofs of a theorem concerning algebraic 
space curve 
\inbook Proceeding of the Eighth National Mathematics Conference (Arya-Mehr Univ. Tech., Teheran, 1977)
\pages 350--370
\yr
\endref

\ref \key H \by R.Hartshorne
\inbook Algebraic Geometry (GTM 52) 
\eds Springer-Verlag 
\publaddr Berlin, 1993
\endref

\ref \key H1 \bysame
\paper Stable vector bundles of rank 2 on $\Bbb P^3$
\jour Math.Ann. \vol 238
\yr 1978
\pages 229--280
\endref

\ref \key HM \by G.Horrocks-D.Mumford
\paper A rank $2$ vector bundle on $\Bbb P^4$ with $15,000$ symmetries
\jour Topology
\vol 12
\yr 1973
\pages 63--81
\endref

\ref \key OSS 
\by C.Okonek-M.Schneider-H.Spindler 
\inbook  Vector Bundles on Complex
Projective Spaces (Prog. Math. 3) 
\eds Bi\-rkh\"auser 
\publaddr Boston, Mass., 1980
\endref

\ref \key O 
\by G.Ottaviani
\paper Some extension of Horrocks criterion to vector bundles
on Grassmannians and quadrics
\jour Ann. Mat. Pura Appl. (4)
\vol 155 
\yr 1989
\pages 317--341
\endref
\endRefs
\enddocument